\documentclass[s5,11pt]{article}
\usepackage{amsmath}
\usepackage{ latexsym}
\newtheorem{theorem}{Theorem}

\newtheorem {definition}{Definition}

\author{ Tord Sj\"odin\footnote{Department of Mathematics and mathematical statistics, Ume\aa\ university, 90187 Ume\aa, Sweden; e-mail, tord.sjodin@umu.se;
phone, +46 907868822}
}
\title{ On Mixtures of Gamma Distributions, Distributions with Hyperbolically  Monotone Densities and Generalized Gamma Convolutions (GGC) }

\begin{document} 
 \maketitle

\begin{abstract} Let $Y$ be a standard Gamma(k) distributed random variable (rv), $k>0$, and let $X$ be an independent positive rv. If $X$ has a hyperbolically monotone density of order $k$ ($HM_k$), then $Y\cdot X$ and $Y/X$ are generalized gamma convolutions (GGC). This extends work by Roynette et al. and Behme and Bondesson. The same conclusions hold when $Y$ is replaced by a finite sum of independent gamma variables with sum of shape parameters at most $k$. Both results are applied to subclasses of GGC.\end{abstract}
\paragraph{ \it AMS 2010 Subject Classsification:} Primary 60E10
Secondary 62E15
\section{Introduction}We consider classes of distributions on the positive real axis. Among these are gamma distributions, distributions with hyperbolically monotone densities ($HM_k$), distributions with hyperbolically completely monotone densities (HCM) and generalized gamma convolutions (GGC). Generalized gamma convolutions were introduced by O. Thorin \cite{T}, \cite{T2} in his study of infinite divisibility of the log-normal distribution. The class GGC consists of limit distributions of convolutions of independent gamma distributions and is closed with respect to (wrt) weak limits and sums and products of independent rvs. A comprehensive study of GGC and its relation to HCM is found in \cite{B1} and \cite{B4}. For more on Bernstein functions and infinite divisibility, see \cite{SSV} and \cite{SH}.\\[1em] 
Recently, Behme and Bondesson \cite{BB} generalized a result of Roynette et al. \cite{RVY} on products $Y\cdot X$ and quotients $Y/X$ of a gamma distributed variable $Y$ and an independent rv $X$ with hyperbolically monotone density. More exactly, they proved that if $0<r\leq k$, $k$ is a positive integer, $Y\sim Gamma(r)$ and $X\sim HM_k$, then $Y\cdot X\sim GGC$ and $Y/X\sim GGC$, \cite{BB}, Theorm 1 and Corollary 1, while Roynette et al.  \cite{RVY} proved the case $k=1$. It was conjectured in \cite{BB} that this result remains true for all $k>0$. We confirm this conjecture and prove that if $Y\sim Gamma(k)$ and $X\sim HM_k$ are independent rvs, then $Y\cdot X\sim GGC$ and $Y/X\sim GGC$, for all $k>0$ (Theorem 1). In addition, we prove the same result when $Y$ is replaced by a finite sum of independent gamma variables with the sum of the shape parameters at most $k$ (Theorem 2). The proofs are based on Bondesson's characterization of $GGC$ and explicit integral expressions for the relevant functions. The results are applied to subclasses of $GGC$ (Theorems 3 and 4). \\[1em]
 The plan of this paper is as follows. Section 2 begins with the standard notation and definitions used in this field and reviews a part of the set up in \cite{BB}, Section 3. We describe our approach in Section 3 and state and prove our main result (Theorem 1). The last two sections contain applications and comments.
 \section{Preliminaries}  This section gives the background to our results and defines the concepts that are needed to state and prove our theorems. A more complete presentation of this field and a survey of many of its results can be found in \cite{B1} and \cite{BB}. A function $f(x)$ is called hyperbolically monotone ($HM_1$) if, for every fixed $u>0$, $h(w)=f(uv)\cdot f(u/v)$ is non-incerasing as a function of $w=v+v^{-1}$. If $k$ is a positive integer, $f$ is called hyperbolically monotone of order $k$ ($HM_k$) if $(-1)^j\cdot D^jh(w)\geq 0$, for $j=1,2,\dots ,k-1$, and $(-1)^{k-1}\cdot D^{(k-1)}h(w)$ decreases. Finally, we say that $f$ is hyperbolically completely monotone if this holds for all positive integers $k$ or, equivalently, that $h(w)$ is completely monotone (CM). See \cite{B1}, Ch. 5 and \cite{BB} for more on these classes. It is proved in \cite{BB} that a density function $f(x)$ is $HM_k$ if and only if, for every $u>0$,  
\begin{equation}h(w)=f(uv)\cdot f(u/v)=\int _{(w,\infty )} (\lambda -w)^{k-1}H_u(d\lambda ),\end{equation}
for $w=v+v^{-1}$ and some nonnegative measure $H_u(d\lambda )$ depending on $u$. We take (1) as our definition of $HM_k$, $k>0$.
\begin{definition} A rv $X \sim HM_k$, if $X$ has density $f(x)$ satisfying (1).\end{definition}
We let $Gamma(k ,r)$ denote the standard class of gamma distributions with density $f(x)= r^k \cdot \Gamma (k )^{-1}\cdot x^{k -1}\cdot e^{-rx}$, $x>0$, $Gamma(k )=Gamma(k ,1)$. A generalized gamma convolution (GGC) is a probability distribution $F$ on $[0,\infty)$ with Laplace transform
$$\phi (s)=\int e^{-sx}\, F(dx)=\exp \big( -as+\int \log (\frac {t}{t+s})\, U(dt)\big),\, s\geq 0,$$
where $a\geq 0$ (called the left extremity) and $U(dt)$ is a nonnegative measure on $(0,\infty )$, with finite mass on all compact subintervals, such that $\int _0^1|\log t| \, U(dt)<\infty $ and $\int _1^\infty t^{-1}\, U(dt)<\infty$, \cite{B1} Ch. 3.
 \\[1em]In the following we prepare for the proof of Theorem 1 and follow the presentation in \cite{BB}, Section 3. Let $X\sim HM_k$, with density $f(x)$ and $Y\sim Gamma(k,r)$ be independent rvs. Then the Laplace transform $\phi (s)$ of the quotient  $Y/X$ is given by
$$\phi (s)=E[e^{-sY/X}] = \int\limits _0^\infty E[e^{-sY/x}]\cdot f(x)\, dx=\int\limits _0^\infty  \bigg(\frac{rx}{rx+s}\bigg)^{k }\cdot f(x)\, dx$$ 
and we get
$$\phi (st)(\phi (s/t)=\int\limits _0^\infty\int\limits _0^\infty  \bigg(\frac{r^2xy}{(rx+st)(ry+s/t)}\bigg)^{k }\cdot f(x)f(y)\, dxdy $$
$$=\int\limits _0^\infty\int\limits _0^\infty \frac{2u}{v} \bigg(\frac{r^2u^2}{(ruv+st)(ru/v+s/t)}\bigg)^{k }\cdot f(uv)f(u/v)\, dudv .$$
Inserting the definition of $HM_k$, putting $a=(ru)/s$ and changing the order of integration, we find that a sufficient condition for $Y/X\sim GGC$ is that the inner most integral
 \begin{equation}J_k=\int\limits _{1/b}^b  \bigg(\frac{ (b-v)(v-b^{-1}) }{(v+ta)(v+t/a) }\bigg)^k\cdot \frac{t^k}{(b-v)(v-b^{-1})}\, dv\end{equation}
   is CM wrt $T=t+\frac{1}{t}$, for all $a>1$ and $b>1$. From now on we assume that $a$ and $b$ are fixed but arbitrary such numbers. The substitution $v\rightarrow v^{-1}$ confirms that the map $t\rightarrow \frac{1}{t}$ leaves $J_k$ invariant and that $J_k$ is a function of $T$. \\[1em]
   The integrand in (2) is a rational function when $k$ is a positive integer and  $J_k $ can, at least in principle, be calculated explicitly. A  straight forward calculation gives 
   $J_1 =(a-a^{-1})^{-1}\cdot \log \big((T+A)/(T+B)\big),$
   where $A=ab+(ab)^{-1}$ and $B=ab^{-1}+a^{-1}b$, and it is easily seen that $J_1 $ is CM wrt $T$, since  $A>B$. In the general integer case, $J_k $ is expressed in terms of certain polynomials $P_k$ and $Q_k$ and the logarithm above, see \cite{BB}, Theorem 1. Our proof for $k>0$ is based on an explicit evaluation of the integral $J_k$ and follows quite different lines.
  \section{Main result}
In this section we state and prove our first result that products and quotients of independent rvs $Y\sim Gamma (k)$ and $X\sim HM_k$, $k>0$, have GGC distributions. This result generalizes, respectively, the cases $k=1$ in \cite{RVY} and $k$ a positive integer in \cite{BB}, to arbitrary $k>0$.b
\begin{theorem} Let $0<k\leq l$ and let $Y\sim Gamma (k)$ and $X\sim HM_l$ be independent rvs. Then $Y\cdot X\sim GGC$ and $Y/X\sim GGC.$
\end{theorem}
 It is sufficient to prove Theorem 1 for $k=l$, since $HM_l\subseteq HM_k$, for $0<k\leq l$, see \cite{BB}. \\[0.5em]
 {\it Proof .} Since $X\sim  HM_k$ if and only if $1/X\sim HM_k$ it is enough to prove that $Y/X\sim HM_k$.  We evaluate $J_k$ by a series of substitutions and first let $$v=(bx+b^{-1})(1+x)^{-1}, \quad dv/dx=(b-b^{-1})\cdot (1+x)^{-2 }$$ and $(b-v)(v-b^{-1})=(b-b^{-1})^2\cdot x\cdot (1+x)^{-2}$. Then $J_k=\int\limits _0^\infty E \, \frac{dx}{x}$, where
 $$E=\bigg(\frac{(b-b^{-1})^2}{(b^2x+\frac{1}{b^2x})\cdot \frac{1}{t}+2T+(x+\frac{t}{x})\cdot t+(bx+\frac{1}{bx}+\beta )\cdot \alpha  }\bigg) ^k\cdot (b-b^{-1})^{-1},
 $$
 where $\alpha =a+a^{-1}$ and $\beta =b+b^{-1}$.  
 The denominator $D$ in $E$ can be rewritten as
 $$D=x\cdot \big( b^2/t+t+\alpha b\big)+\frac{1}{x}\cdot \big(  1/(b^2t)+t+\alpha /b\big)+2T+\alpha\beta  .
 $$
 and a standard substitution, putting the second term in $D$ equal to $1/\rho$, transforms $D$ into
 $$D_1 =\rho^{-1} +\rho\cdot (T+A)(T+B)+2T+\alpha\beta  ,$$
 where $A=ab+(ab)^{-1}$, $B=ab^{-1}+a^{-1}b$.
   Now let $\delta =\rho \cdot \sqrt{(T+A)(T+B)}$ and we get
 \begin{equation} J_k=\int\limits _0^\infty \bigg(\frac{(b-b^{-1})^2}{(\delta +\frac{1}{\delta})
  \sqrt{(T+A)(T+B)}+2T+\alpha \beta  }\bigg)^k\cdot (b-b^{-1})^{-1}\, \frac{d\delta}{\delta},
 \end{equation}
 which is our final expression for $J_k$.
 The denominator in (3) is now easily seen to have a CM derivative and hence $J_k$ is CM wrt $T$, see \cite{F}, p. 441. The proof of Theorem 1 is complete.\hfill $\Box$ 
\section{Sums of independent gamma variables} In this section we replace the gamma variable $Y$ in Theorem 1 by a finite sum of independent gamma variables with sum of shape parameters at most $k$, as was suggested in \cite{BB}. We prove that the conclusions of Theorem 1 remain true (Theorem 2) and apply the result to subclasses of GGC in the next section.
\begin{theorem}Let $Y=Y_1+Y_2+\cdots +Y_n$ be a finite sum of
 independent gamma variables, $Y_i\sim Gamma (k_i,r_i)$, $1\leq i\leq n$, such that $k_1+k_2+\cdots +k_n\leq k$ and let $X\sim HM_k$ be an independent rv. Then $Y\cdot X\sim GGC$ and $Y/X\sim GGC.$
\end{theorem}
{\it Proof.} It is no loss of generality to assume that $k_1+k_2+\cdots +k_n= k$. Let $Y$ and $X$ be as in the theorem, then as above
$$\phi (s)=E[e^{-sY/X}] = \int\limits _0^\infty E[e^{-sY/x}]\cdot f(x)\, dx=\int\limits _0^\infty\prod\limits_{i=1}^n \bigg(\frac{r_ix}{r_ix+s}\bigg)^{k_i}\cdot f(x)\, dx$$
and 
$$\phi (st)(\phi (s/t)=\int\limits _0^\infty\int\limits _0^\infty\frac{2u}{v}\prod\limits_{i=1}^n \bigg(\frac{b_i^2u^2}{(b_iuv+st)(b_iu/v+s/t)}\bigg)^{k_i}\cdot f(uv)f(u/v)\, dudv.$$
 We now proceed  as in the proof of Theorem 1 and find that it is sufficient 
  to prove that the following analogue of (2) $$J_k=\int\limits_{1/b}^b\prod\limits_{i=1}^n\bigg(\frac{(b-v)(v-1/b)}{(v+ta_i)(v+t/a_i)}\bigg)^{k_i}\cdot \frac{t^k}{(b-v)(v-1/b)}\, dv$$
is CM wrt $T=t+t^{-1}$, where $a_i=(r_iu)/s$, $1\leq i\leq n$. Then the substitution $v=(bx+b^{-1})(1+x)^{-1}$ gives
$$J_k= \int\limits_0^\infty \prod\limits_{i=1}^n\bigg(\frac{(b-b^{-1})^2}{x\cdot (\frac{b^2}{t}+t+\alpha_ib)+\frac{1}{x}\cdot (\frac{1}{b^2t}+t+\frac{\alpha_i}{b})+2T+\alpha_i\beta}\bigg)^{k_i}\cdot \frac{1}{b-b^{-1}} \frac{dx}{x}.$$
Denote the the i-th factor in the denominator of the integrand by $S_i$, then the integrand is CM wrt the variables $S_i$, $1\leq i\leq n$, and by Bernstein's theorem there is a nonnegative measure $\nu (d\lambda )$ such that $J_k=\int\limits_0^\infty \frac{dx}{x}\int\limits _0^\infty \nu (d\lambda )\, e^{-E}$, where $E= \sum \limits_1^n \lambda_iS_i$. Next we change the order of integration in $J_k$, denote the inner integral by $I_k$ and set out to prove that $I_k$ is CM wrt $T$, for all choices of $\lambda_i$, $1\leq \lambda_i\leq n$. We assume, without loss of generality, that $\sum\limits _1^n \lambda_i=1$, define $\alpha =\sum\limits _1^n \lambda_i\alpha_i $ and choose $a>1$ such that $\alpha =a+a^{-1}$. Then
$$ E=x\cdot (b^2/t+t+\alpha b)+x^{-1}\cdot (1/(b^2t)+t+ \alpha /b)+2T+\alpha\beta $$
   and performing the same substitutions as in the proof of Theorem 1 gives
 $$I_k=\int\limits_0^\infty \exp \bigg(-\big((\delta +\delta^{-1})
  \sqrt{(T+A)(T+B)}+2T+\alpha \beta \big)\bigg)   \frac{d\delta}{\delta},
  $$
  with $\alpha $ and $ a$ defined above and $A=ab+(ab)^{-1}$, $B=ab^{-1}+a^{-1}b$. This shows that $J_k$ is CM wrt $T$ and completes the proof of Theorem 2.\hfill $\Box$
 \section{Applications}Let $\mathcal A$ and $\mathcal B$ be classes of distributions and denote by $\mathcal A\times\mathcal B$ the class of products $Y\cdot X$ for independent rvs $Y\sim \mathcal A$ and $X\sim \mathcal B$ and similarily for quotients. Then  $Gamma (k)\times HM_k\subset GGC$ and $Gamma (k)/ HM_k\subset GGC$, for $k>0$, by Theorem 1. If $Z\sim GGC$ is independent of $Y$and $X$, then $(Z\cdot  Y)\cdot X=Z\cdot ( Y\cdot X)\sim GGC$ and $(Z\cdot  Y)/ X=Z\cdot ( Y/ X)\sim GGC$, since $GGC $ is closed wrt independent products, \cite{B4},   Theorem 1, and we have the following result.\begin{theorem}  For any $k>0$,  $(GGC\times Gamma (k))\times HM_k\subset GGC$ and $(GGC\times Gamma (k))/HM_k\subset GGC.$
 \end{theorem} If $\mathcal H_k$ is the largest class of probability distributions such that $ Gamma(k)\times \mathcal H_k\subset GGC$, then $HM_k\times GGC \subset \mathcal H_k$. Similarily, if $\mathcal G_k$ is the largest class of probability distributions such that $\mathcal G_k\times HM_k\subset GGC$, then $GGC\times Gamma(k)\subset \mathcal G_k$. We can use Theorem 2 to improve on this. Let $GGC(k)$ be the class of rvs $X\sim GGC$ with left extremity zero and total $U-$measure at most $k$.
 \begin{theorem} For any $k>0$, $GGC(k)\subset \mathcal G_k$.
 \end{theorem}
 {\it Proof.} Every $X\sim GGC(k)$ is the weak limit of finite sums of independent gamma variables with $U-$measure at most $k$ and the conclusion follows from Theorem 2.\hfill $\Box$
      \section{Final comments} The extension of the integer case of Theorem 1 in \cite{BB} to hold for all $k\geq 1$, and even for all $k>0$, was suggested to me by Professor Bondesson. Our proof shows that, in contrast to what is believed in \cite{BB}, the function $J_k$ (3) can be explicitely calculated also in the non-integer case. The truth of Theorem 2 and the possibility that $GGC(k)$ is contained in $\mathcal G_k$ is mentioned in \cite{BB}.\section{Acknowledgements} The author thanks Professor Lennart Bondesson for introducing me to the problems studied here and for valuable comments and discussions.

 \end{document}